# Bifurcating Continued Fractions

**Ashok Kumar Gupta**
**Department of Electronics and Communication**
**Allahabad University, Allahabad - 211 002, INDIA**
**(Email address: akgjkiapt@hotmail.com)**

**Ashok Kumar Mittal**
**Department of Physics**
**Allahabad University, Allahabad – 211 002, India**
**(Email address: mittal_a@vsnl.com)**

**Abstract:**

The notion of 'bifurcating continued fractions' is introduced. Two coupled sequences of non-negative integers are obtained from an ordered pair of positive real numbers in a manner that generalizes the notion of continued fractions. These sequences enable simple representations of roots of cubic equations. In particular, remarkably simple and elegant 'bifurcating continued fraction' representations of Tribonacci and Moore numbers, the cubic variations of the 'golden mean', are obtained. This is further generalized to associate $m$ non-negative integer sequences with a set of $m$ given real numbers so as to provide simple 'bifurcating continued fraction' representation of roots of polynomial equations of degree $m$+1.

## 1. Introduction

Continued fractions provide much insight into mathematical problems, particularly into the nature of numbers.[1] In the computer field, continued fractions are used to give approximations to various complicated functions, and once coded for the electronic machines, give rapid numerical results valuable to scientists. A reference to continued fractions is found in the works of the Indian mathematician Aryabhatta. John Wallis used for the first time the name “continued fraction” in his book *Arithmetica Infinitorium*, published in 1655. Christian Huygens used continued fractions for the purpose of approximating the correct design for the toothed wheels of a planetarium. Euler, Lambert and Lagrange were prominent amongst those who developed the theory of continued fractions.

Any eventually periodic continued fraction represents a quadratic irrational. Conversely, Lagrange's theorem asserts that the continued fraction expansion of every quadratic irrational is eventually periodic. A purely periodic continued fraction represents a quadratic irrational of a special kind called a reduced quadratic irrational. A quadratic irrational is said to be reduced if it is greater than 1 and the other root of the quadratic equation that it satisfies, lies between -1 and 0. Conversely, the continued fraction expansion of a reduced quadratic irrational is purely periodic.

The continued fraction expansion consisting of the number 1 repeated indefinitely represents the 'golden mean'. This satisfies the quadratic equation $x^2 = x + 1$. The convergents of the continued fraction are obtained as the ratio of the successive terms of the Fibonacci sequence. Each term of the Fibonacci sequence is obtained by summing the previous two terms of the sequence.

A straightforward generalization of the Fibonacci sequence is that of the Tribonacci sequence, in which each term is obtained by summing the previous three terms of the sequence. Tribonacci number, the limiting ratio of the successive terms of the Tribonacci sequence satisfy the Tribonacci equation $x^3 = x^2 + x + 1$. However, the conventional continued fraction expansion of the Tribonacci number offers no satisfying pattern that would reflect the simple generalization from the 'golden mean'.

In this paper, a method to generalize the continued fraction expansion is suggested so as to remove this deficiency of the conventional continued fractions. This method associates two coupled non-negative integer sequences with two real numbers and based on these sequences a 'bifurcating continued fraction' expansion is obtained. The 'bifurcating continued fraction' resembles the conventional continued fraction, except that both the numerator and the denominator bifurcate in a “Fibonacci Tree” like manner, whereas in conventional continued fractions only the denominator bifurcates. In particular, two sequences consisting of the number 1 repeated indefinitely represent the Tribonacci number. Other cubic variants[2] of the 'golden mean', like the Moore number satisfying the cubic equation $x^3 = x^2 + 1$, also find very simple and elegant representation as a 'bifurcating continued fraction'. Bifurcating continued fraction representations reveal the secret beauty of many numbers.

This is further generalized to associate *m* non-negative integer sequences with a set of *m* given real numbers so as to provide simple 'bifurcating continued fraction' representation of roots of polynomial equations of degree $m + 1$.

## 2. Generalizing the Continued Fractions

Given a positive real number $\alpha$, the continued fraction expansion of $\alpha$ is given by a sequence of non-negative integers $[a_0, a_1, \ldots\ldots a_i \ldots.]$, obtained by the recurrence relation

$$a_i = \text{int}(\alpha_i),$$
$$\alpha_{i+1} = 1/(\alpha_i - a_i)\,, \qquad (1)$$

where $\alpha_0 = \alpha$.

Our generalization of the method of continued fraction is based on a generalization of equation (1). In this generalization we obtain an ordered pair of integer sequences $\{\{a_0, a_1,\ldots,a_i, \ldots\}, \{b_0, b_1,\ldots.,b_i,\ldots\}\}$ from a given ordered pair $\{\alpha,\beta\}$ of positive real numbers by using the recurrence relations

$$a_i = \text{int}(\alpha_i),$$
$$b_i = \text{int}(\beta_i),$$
$$\alpha_{i+1} = 1/(\beta_i - b_i),$$
$$\beta_{i+1} = (\alpha_i - a_i)/\beta_i - b_i), \quad (2)$$

where $\alpha_0 = \alpha$ and $\beta_0 = \beta$.

As an example, consider the pair of numbers $\alpha = 2^{1/3}$ and $\beta = 2^{2/3}$. Given this pair of numbers the integer sequences obtained from equation (2) are given by

$\{a_i\}$ = { 1 1 2 1 2 1 2 1 2 1 2 1 2 1 2 1 2 1 2 1 ………}

= {1 (1 2) }

$\{b_i\}$ = { 1 0 1 0 1 0 1 0 1 0 1 0 1 0 1 0 1 0 1 0 ……..}

= { (1 0) }

These sequences have a period 2. The simplicity of these sequences is in sharp contrast to the conventional continued fraction representations of these numbers given below:

$2^{1/3}$ = [1, 3, 1, 5, 1, 1, 4, 1, 1, 8, 1, 14, 1, 10, 2, 1, 4, 12, 2, 3, 2, 1, 3, 4, 1, 1, 2, 14, 3, 12, 1, 15, 3, 1, 4, 534, 1, 1, 5, 1, 1, 121, 1, 2, 2, 4, 10, 3, 2, 2, 41, 1, 1, 1, 3, 7, 2, 2, 9, 4, 1, 3, 7, 6, 1, 1, 2, 2, 9, 3, 1, 1, 69, 4, 4, 5, 12, 1, 1, 5, 15, 1, 4, 1, 1, 1, 1, 1, 89, 1, 22, 186, 6, 2, 3, 1, 3, 2, 1, 1, 5, 1, 3, 1, 8, 9, 1, 26, 1, 7, 1, 18, 6, 1, 372, 3, 13, 1, 1, 14, 2, 2, 2, 1, 1, 4, 3, 2, 2, 1, 1, 9, 1, 6, 1, 38, 1, 2, 25, 1, 4, 2, 44, 1, 22, 2, 12, 11, 1, 1, 49, 2, 6, 8, 2, 3, 2, 1, 3, 5, 1, 1, 1, 3, 1, 2, 1, 2, 4, 1, 1, 3,…..]

$2^{2/3}$ = [1, 1, 1, 2, 2, 1, 3, 2, 3, 1, 3, 1, 30, 1, 4, 1, 2, 9, 6, 4, 1, 1, 2, 7, 2, 3, 2, 1, 6, 1, 1, 1, 25, 1, 7, 7, 1, 1, 1, 1, 266, 1, 3, 2, 1, 3, 60, 1, 5, 1, 8, 5, 6, 1, 4, 20, 1, 4, 1, 1, 14, 1, 4, 4, 1, 1, 1, 1, 7, 3, 1, 1, 2, 1, 3, 1, 4, 4, 1, 1, 1, 3, 1, 34, 8, 2, 10, 6, 3, 1, 2, 31, 1, 1, 1, 4, 3, 44, 1, 45, 93, 12, 1, 7, 1, 1, 5, 12, 1, 1, 2, 4, 19, 1, 12, 1, 16, 1, 8, 1, 1, 2, 1, 745, 1, 1, 1, 6, 3, 1, 6, 1, 2, 2, 2, 3, 2, 6, 1, 5, 20, 1, 2, 1, 78, 2, 1, 12, 2, 2, 4, 22, 2, 11, 4, 6, 23, 99, 1, 12, 4, 4, 1, 1, 2, 7, 2, 1, 4, 1, 1, 2, …]

The numbers $\alpha$ and $\beta$ can be recovered from the two integer sequences $\{a_i\}$ and $\{b_i\}$ by using the recurrence relations

$$\alpha_i = a_i + (\beta_{i+1} / \alpha_{i+1}),$$
$$\beta_i = b_i + (1/\alpha_{i+1}) \quad (3)$$

Let us introduce the following notation:

$$\alpha_n = [\, a_n, b_n; a_{n+1}, b_{n+1}; a_{n+2}, b_{n+2}; \ \ldots\ldots\,]$$

and

$$\beta_n = [\, b_n; a_{n+1}, b_{n+1}; a_{n+2}, b_{n+2}; a_{n+3}, \ \ldots\ldots\,]$$

Then equations (3) can be expressed by

$$[\, a_i, b_i; a_{i+1}, b_{i+1}; a_{i+2}, b_{i+2};\ \ldots\ldots] = a_i + \frac{[b_{i+1}; a_{i+2}, b_{i+2}; a_{i+3}, b_{i+3};\ \ldots\ldots]}{[a_{i+1}, b_{i+1}; a_{i+2,} b_{i+2}; a_{i+3}, b_{i+3};\ldots.]}$$

$$[\, b_i; a_{i+1}, b_{i+1}; a_{i+2}, b_{i+2};\ \ldots\ldots] = b_i + \frac{1}{[a_{i+1}, b_{i+1}; a_{i+2,} b_{i+2}; a_{i+3}, b_{i+3};\ldots.]}$$

Then,

$$\alpha = \alpha_0 = a_0 + \frac{[b_1; a_2, b_2; a_3, b_3;\ \ldots\ldots]}{[a_1, b_1; a_2, b_2; a_3, b_3;\ldots.]}$$

$$= a_0 + \cfrac{b_1 + \cfrac{1}{[a_2, b_2; a_3, b_3;\ \ldots\ldots\ldots]}}{a_1 + \cfrac{[b_2; a_3, b_3;\ \ldots\ldots\ldots\ldots]}{[a_2, b_2; a_3, b_3;\ \ldots\ldots\ldots]}}$$

$$= a_0 + \cfrac{b_1 + \cfrac{1}{a_2 + \cfrac{b_3 + \cfrac{1}{[\, a_4, b_4; a_5, b_5;\ldots\ldots]}}{a_3 + \cfrac{[\, b_4; a_5, b_5;\ldots\ \ldots\ldots]}{[\, a_4, b_4; a_5, b_5;\ldots\ldots]}}}}{a_1 + \cfrac{b_2 + \cfrac{1}{a_3 + \cfrac{[\, b_4; a_5, b_5;\ldots\ \ldots\ldots]}{[\, a_4, b_4; a_5, b_5;\ldots\ldots\ldots]}}}{a_2 + \cfrac{b_3 + \cfrac{1}{[\, a_4, b_4; a_5, b_5;\ldots\ldots]}}{a_3 + \cfrac{[\, b_4; a_5, b_5;\ldots\ \ldots\ldots]}{[\, a_4, b_4; a_5, b_5;\ldots\ldots]}}}}$$

Given the integer sequences $\{a_i\}$ and $\{b_i\}$, one can write down the bifurcating continued fractions by the bifurcating rule:

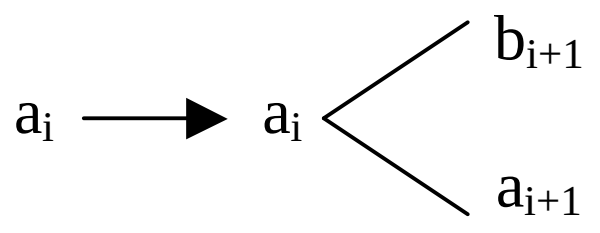

and

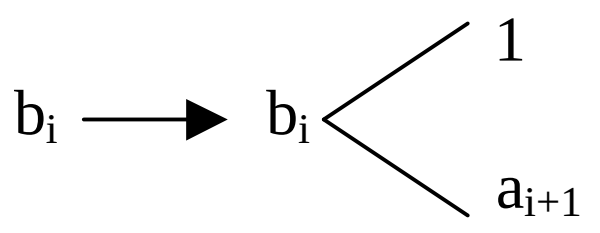

where

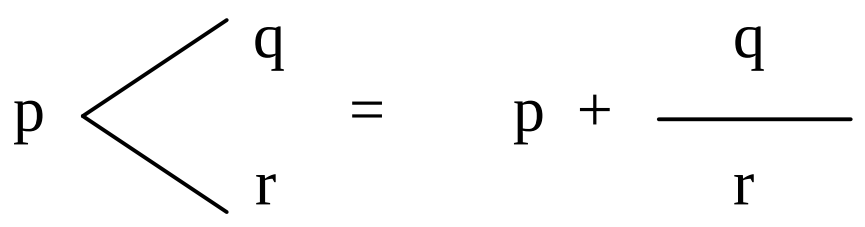

To get $\alpha$, one begins with $a_0$ and to get $\beta$ one begins with $b_0$. In this notation

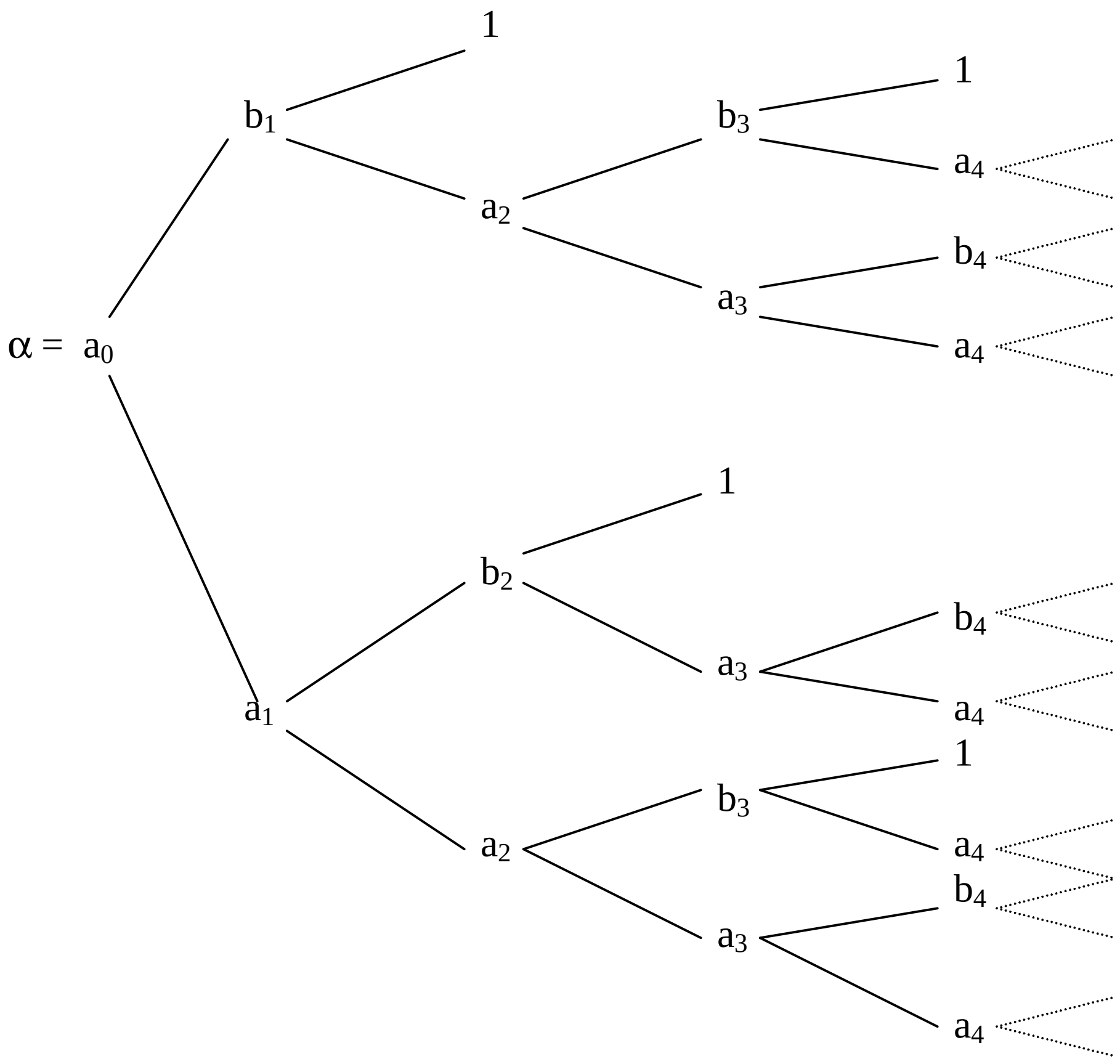

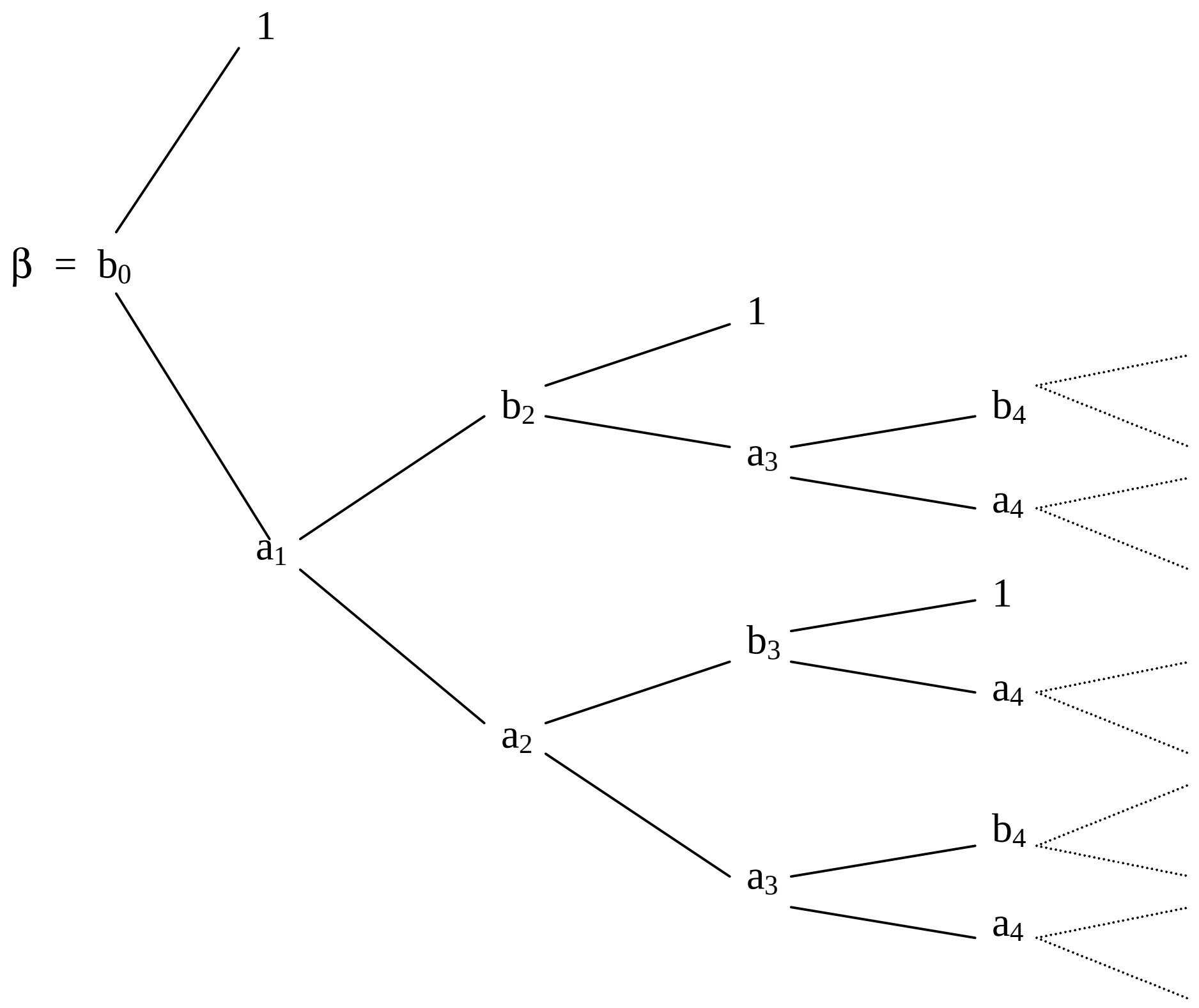

The Fibonacci tree structure of these expansions is clear.

## 3. Period One Bifurcating Continued Fractions

The simplest bifurcating continued fractions are those with period one. In this case the sequence $\{a_i\}$ consists of an integer 'a' repeated indefinitely and the sequence {b} consists of integer 'b' repeated indefinitely. Let $\alpha$ and $\beta$ be the pair of numbers represented by this pair of constant sequences. Then it follows from the equations (3) that $\alpha$ and $\beta$ satisfy the equations

$$\alpha = a + (\beta/\alpha)$$

and (4)

$$\beta = b + (1/\alpha)$$

Thus $\alpha$ and $\beta$ satisfy the cubic equations

$$\alpha^3 = a\,\alpha^2 + b\,\alpha + 1 \qquad (5)$$

and

$$\beta^3 = 2b\,\beta^2 - (a + b^2)\,\beta + (ab + 1) \qquad (6)$$

We represent them in the Fibonacci tree notation as follows:

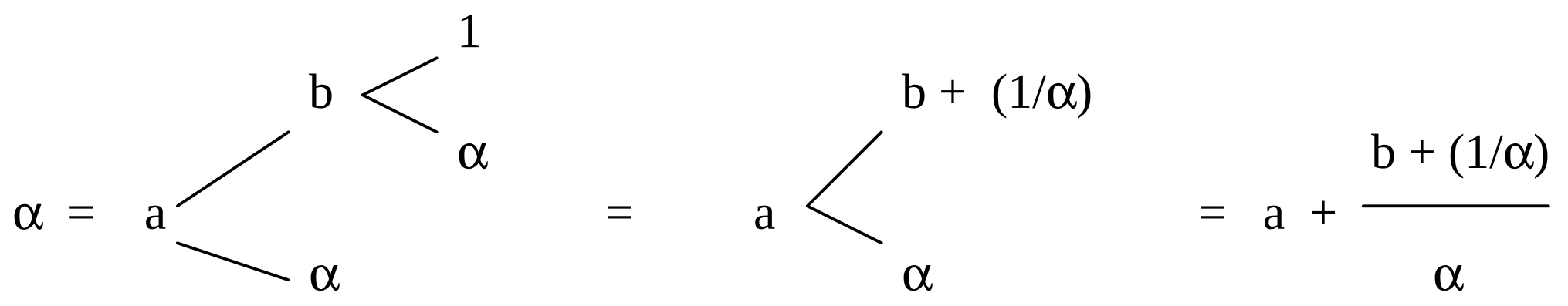

Therefore

$$\alpha^3 = a\,\alpha^2 + b\,\alpha + 1.$$

and

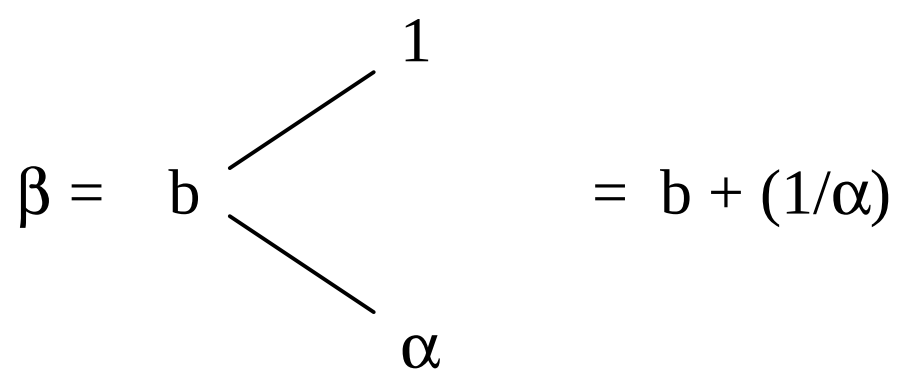

For example, we consider the two simplest cases.

(i) If we take a = b = 1, $\alpha$ is the Tribonacci number. The Tribonacci sequence defined by

$t_0 = 0;\ t_1 = 0;\ t_2 = 1$

$t_{n+2} = t_n + t_{n+1} + t_{n+2},\ \ n > 0$

is the generalization of the famous Fibonacci sequence. The limiting ratio

$$\alpha = \lim \frac{t_{n+1}}{t_n} \text{ as } n \to \infty$$

is the real solution of $\chi^3 = \chi^2 + \chi + 1$. It is equal to 1.83928675521416….

It's simple continued fraction representation is given by [1, 1, 5, 4, 2, 305, 1, 8, 1, 2, 1, 3, 1, 2, 2, 1, 1,…..]. In sharp contrast with the very simple continued fraction of the 'golden number', the continued fraction of the Tribonacci number exhibits no discernible pattern. In view of the fact that the Tribonacci number is a straightforward generalization of the 'golden number' this is not very satisfactory. Bifurcating continued fractions provide a satisfactory way of representing the Tribonacci number in a manner that generalizes the continued fraction representation of the 'golden mean'. This representation is given below:

$$\alpha = 1 + \cfrac{1 + \cfrac{1 + \cfrac{1}{1 + \cfrac{1}{[\,1, 1;\, 1, 1; \ldots\ldots]}}}{1 + \cfrac{[\,1;\, 1, 1; \ldots\ldots]}{[\,1, 1;\, 1, 1; \ldots\ldots]}}}{1 + \cfrac{1 + \cfrac{1}{1 + \cfrac{[\,1;\, 1, 1; \ldots\ldots]}{[\,1, 1;\, 1, 1; \ldots\ldots]}}}{1 + \cfrac{1 + \cfrac{1}{[\,1, 1;\, 1, 1; \ldots\ldots]}}{1 + \cfrac{[\,1;\, 1, 1; \ldots\ldots]}{[\,1, 1;\, 1, 1; \ldots\ldots]}}}}$$

$$= 1.83928675521416\ldots$$

It is the solution of $\alpha^3 = \alpha^2 + \alpha + 1$. The bifurcating continued fraction of Tribonacci constant given above is a generalization of the simple continued fraction expansion of Golden Constant $\tau = 1.6180339..$

$$\tau = \cfrac{1}{1 + \cfrac{1}{1 + \cfrac{1}{1 + \cfrac{1}{1 + \ldots\ldots}}}}$$

The fifth convergent of $\alpha$ is given as

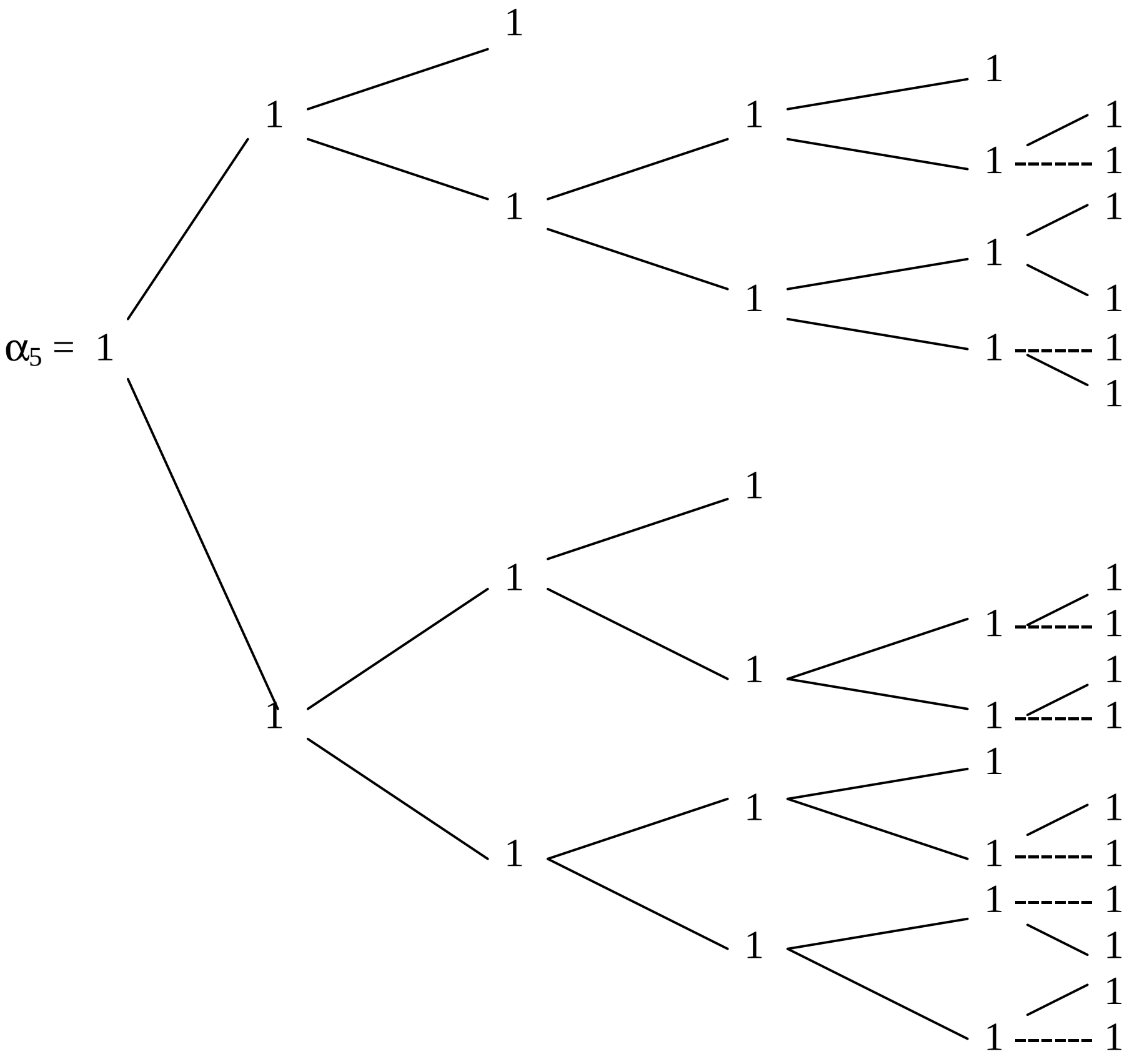

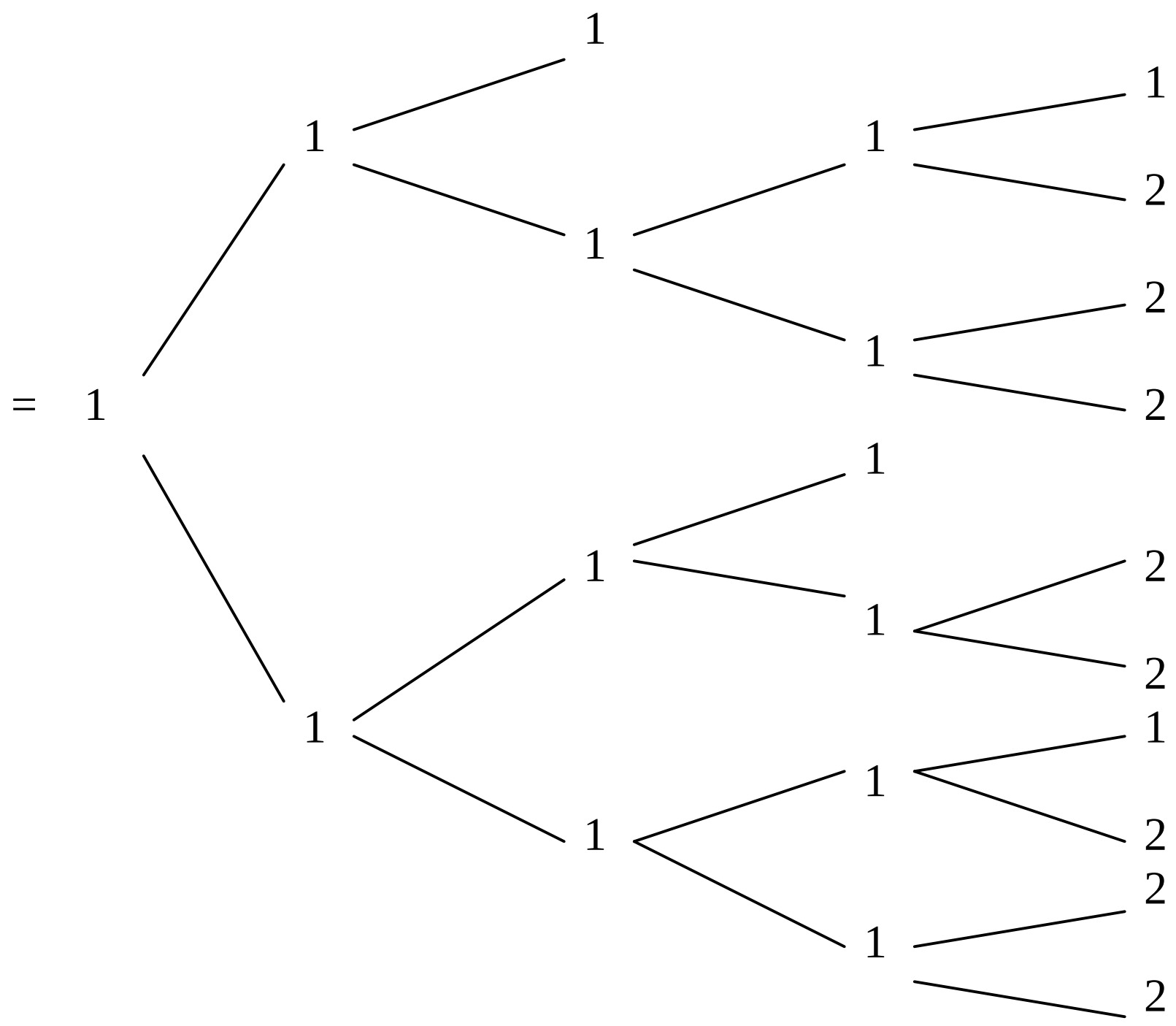

=

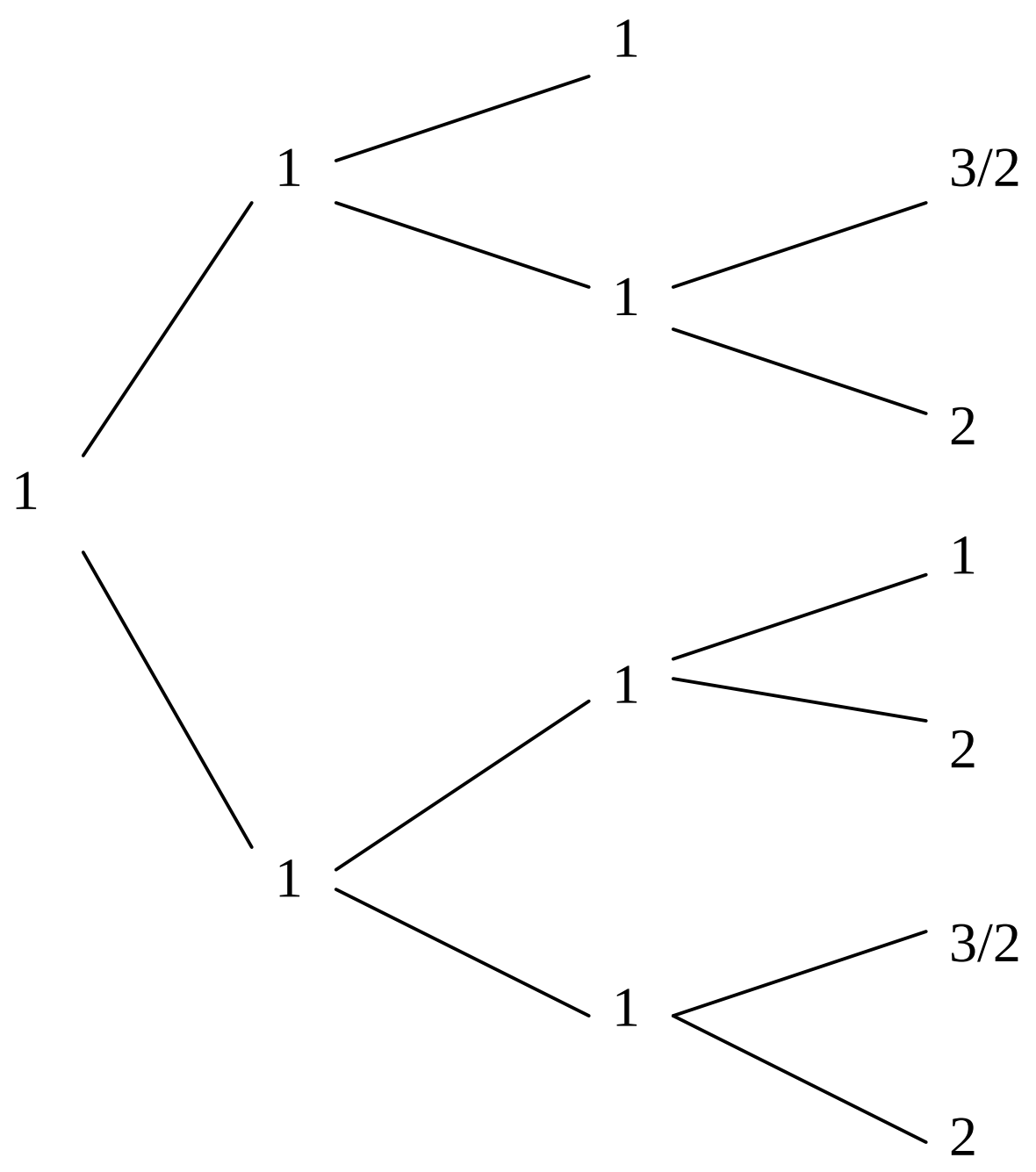

=

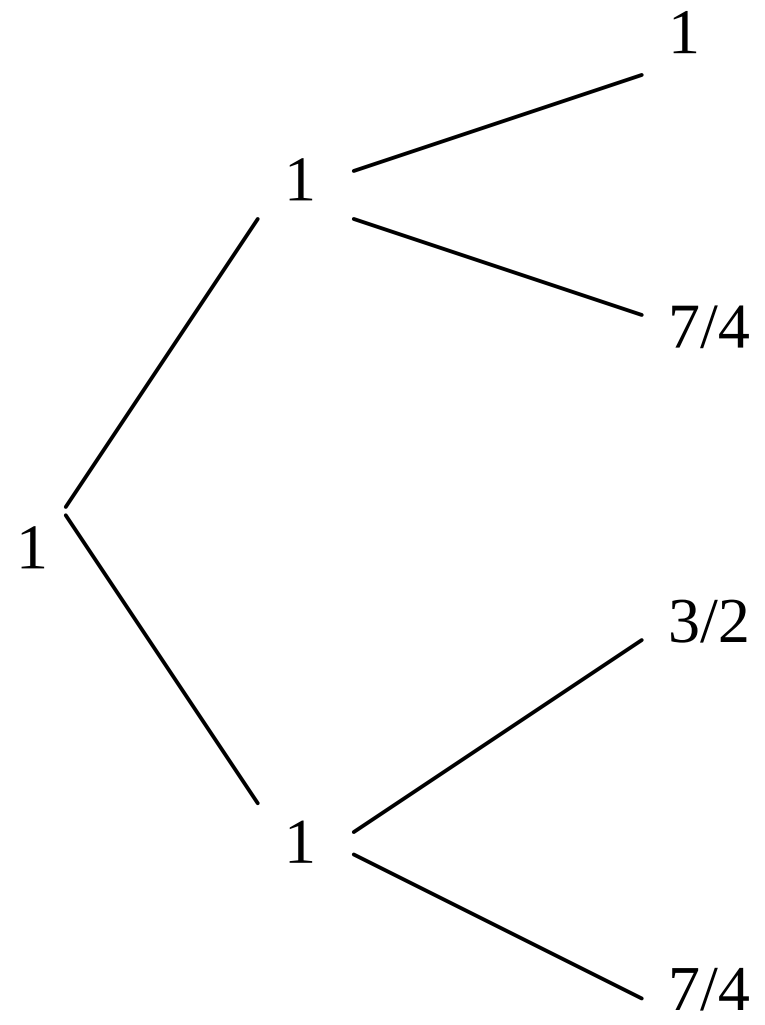

=

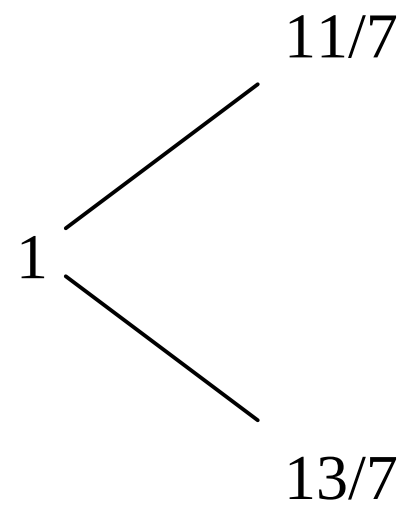

= 24/13 = 1.846

(ii) For a = 1, b = 0, $\alpha$ is the Moore number ( $\alpha$ = 1.4655712318..), which is the solution of the cubic $\chi^3 - \chi^2 - 1 = 0$.

Using the graphical technique it is easy to show that the real numbers represented by BFC of period 2 sequences satisfy a cubic equation. It remains to be seen if Lagrange's theorem can be generalized to cubic irrationals.

## 4. Further Generalization

The generalization to still higher order bifurcating continued fractions (BCF) is straightforward. We can associate to an *m* – tuple $\{\alpha^{(1)}, \alpha^{(2)}, \ldots, \alpha^{(m)}\}$ an m-tuple of integer sequences $\{\{a^{(1)}_0, a^{(1)}_1, \ldots, a^{(1)}_i, \ldots\}, \{a^{(2)}_0, a^{(2)}_1, \ldots, a^{(2)}_i, \ldots\}, \ldots, \{a^{(k)}_0, a^{(k)}_1, \ldots, a^{(k)}_i, \ldots\}, \ldots, \{a^{(m)}_0, a^{(m)}_1, \ldots, a^{(m)}_i, \ldots\}\}$ obtained from the recurrence relations

$$a^{(k)}_i = \mathrm{int}(\alpha^{(k)}_i),$$
$$\alpha^{(1)}_{i+1} = 1/(\alpha^{(m)}_i - a^{(m)}_i),$$
$$\alpha^{(k+1)}_{i+1} = (\alpha^{(k)}_i - a^{(k)}_i)/(\alpha^{(m)}_i - a^{(m)}_i), \qquad k = 1,2, \ldots m-1 \qquad (7)$$

As an example, let $\alpha^{(1)} = 2^{1/4}$, $\alpha^{(2)} = 2^{1/2}$, $\alpha^{(3)} = 2^{3/4}$. Given this 3-tuple of numbers, the integer sequences obtained by application of recurrence relations (7) are

$\{a^{(1)}_i\}$ = { 1 1 1 2 1 1 2 1 1 2 1 1 2 …} = { 1 (1 1 2) }

$\{a^{(2)}_i\}$ = { 1 0 0 1 0 0 1 0 0 1 0 0 ….} = { (1 0 0) }

$\{a^{(3)}_i\}$ = { 1 0 0 1 0 0 1 0 0 1 0 0 ….} = { (1 0 0) }

The numbers $\alpha^{(k)}$ can be recovered from the m-tuple of integer sequences from the recurrence relations:

$$\alpha^{(k)}_i = a^{(k)}_i + (\alpha^{(k+1)}/\alpha^{(1)}_{i+1}), \qquad k = 1,2, \ldots m-1$$
$$\alpha^{(m)}_i = a^{(m)}_i + (1/\alpha^{(1)}_{i+1}), \qquad (8)$$

If we take all unit sequences,

$\{ a_1 \}$ = { 1 1 1 1 1 1 1 1 …..}
$\{ a_2 \}$ = { 1 1 1 1 1 1 1 1 …..}
$\{ a_3 \}$ = { 1 1 1 1 1 1 1 1 …..}

then we get

$\alpha^{(1)}$ = 1.92756197548094..
$\alpha^{(2)}$ = 1.78793319384454..
$\alpha^{(3)}$ = 1.51879006367738..

Here $\alpha^{(1)}$ is the Tetranacci number which is the real root of $x^4 = x^3 + x^2 + x + 1$.

## 5. Acknowledgement

We are grateful to Prof. D. N. Verma, Tata Institute of Fundamental Research, India, for the a suggestion that something like eq. (2) may be useful in the generalization of continued fractions.